\newif\ifspringer
\ifspringer
\documentclass[sn-mathphys]{sn-jnl}
\else
\documentclass[twoside,a4paper]{amsart}
\fi
\usepackage{amssymb}
\usepackage{mathtools}

\PassOptionsToPackage{pdfauthor={Taghreed Alqurashi and Vladimir V. Kisil},%
    pdftitle={Metamorphism as a covariant transform for the SSR group},%
    pdfsubject={mathematics},%
    backref=page,%
  pdfkeywords={symmetries, metamorphism, covariant transform, integral transform}
}{hyperref}
\usepackage{hyperref}
\ifspringer
\providecommand{\amscite}[3]{\cite[#3]{#1}}
\providecommand{\citelist}{}
\providecommand{\cites}[1]{\cite{#1}}
\else
\IfFileExists{eulervm.sty}{%
\usepackage{eulervm}}{}
\usepackage[nobysame]{myamsrefs}
\providecommand{\amscite}[3]{\cite{#1}*{#3}}
\fi

\newtheorem{thm}{Theorem}[section]
\newtheorem{prop}[thm]{Proposition}

\newtheorem{cor}[thm]{Corollary}

\theoremstyle{definition}

    \theoremstyle{remark}
    \newtheorem{rem}[thm]{Remark}

\usepackage{accents}
\newcommand{\metamorph}[1]{\accentset{\smash{\raisebox{-0.12ex}{$\scriptscriptstyle\approx$}}}{#1}\rule{0pt}{2.3ex}}

\providecommand{\half}{{\textstyle\frac{1}{2}}}

   \PassOptionsToPackage{british}{babel}
\providecommand{\norm}[2][\relax]{\left\|#2\right\|\ifx#1\relax\else_{#1}\fi}
\providecommand{\modulus}[2][\relax]{\left| #2 \right|\ifx#1\relax\else_{#1}\fi}
\providecommand{\oper}[1]{\mathcal{#1}}
\providecommand{\algebra}[1]{\ensuremath{\mathfrak{#1}}}

\providecommand{\Space}[3][]{\ifx#2R\ifx#1e \mathbb{C}^{#3} \else
\ifx#1p \mathbb{D}^{#3} \else
\ifx#1h \mathbb{O}^{#3} \else
\ifx#1\sigma \mathbb{A}\!^{#3} \else
\ensuremath{\mathbb{#2}^{#3}_{#1}{}} \fi \fi \fi \fi \else
\ensuremath{\mathbb{#2}^{#3}_{#1}{}} \fi}
\providecommand{\FSpace}[3][]{\ensuremath{\ifx#2l \ell_{#3}^{#1}{}\else
  \mathsf{#2}_{#3}^{#1}{}\fi}}

\providecommand{\uir}[3][0]{\ifcase #1{\rho^{#2}_{#3}}%
\or {\breve{\rho}^{#2}_{#3}}%
\or {\tilde{\rho}^{#2}_{#3}}\fi}
\providecommand{\scalar}[3][\relax]{\left\langle #2, #3
        \right\rangle\ifx#1\relax\else_{#1}\fi}

\providecommand{\SL}[1][2]{\ensuremath{\mathrm{SL}_{#1}(\Space{R}{})}}

\providecommand{\rmi}{\mathrm{i}}
\providecommand{\rme}{\mathrm{e}}
\providecommand{\rmd}{\mathrm{d}}

\providecommand{\myhbar}{\hslash}
\makeatletter
\newcommand{\@hslashslash}{%
  \raisebox{-0.2ex}{%
    \scalebox{1}{%
      \rotatebox[origin=c]{20}{$\mathchar'26$%
      }%
    }%
  }%
}
\newcommand{\@kslash}[2]{%
  {%
   \vphantom{#2}%
   \ooalign{\kern#1em\smash{\@hslashslash}\hidewidth\cr
     $#2$\cr
   }%
   \kern.05em
  }%
}
\newcommand\kslash{\mathchoice
  {\hbox{\@kslash{.01}{k}}}
  {\hbox{\@kslash{.01}{k}}}
  {\hbox{\fontsize{\sf@size}{\sf@size}\selectfont\@kslash{.05}{k}}}
  {\hbox{\fontsize{\ssf@size}{\ssf@size}\selectfont\@kslash{.1}{k}}}
}

\providecommand{\myh}{h}

\begin{document}

\title[Metamorphism as a covariant transform]{Metamorphism as a covariant transform\\ for the SSR group}

\ifspringer
\else

\author[T. Alqurashi and V.V. Kisil]%
{Taghreed Alqurashi and \href{http://www1.maths.leeds.ac.uk/~kisilv/}{Vladimir V. Kisil}}
\thanks{On  leave from Odessa University.}

\address{%
School of Mathematics\\
University of Leeds\\
Leeds LS2\,9JT\\
UK
}
\address{Mathematics Department\\ Albaha University\\ Faculty of Science in Almakhwah\\
Saudi Arabia}
\email{talqorashi@bu.edu.sa}
\email{\href{mailto:V.Kisil@leeds.ac.uk}{V.Kisil@leeds.ac.uk}}

\urladdr{\url{http://www1.maths.leeds.ac.uk/~kisilv/}}

\date{\today}
\begin{abstract}
  Metamorphism is a recently introduced integral transform, which is useful in solving partial differential equations. Basic properties of metamorphism can be verified by direct calculations. In this paper we present metamorphism as a sort of covariant transform and derive its most important features in this way. Our main result is a characterisation of metamorphism's image space. Reading this paper does not require advanced knowledge of group representations or theory of covariant transform.
\end{abstract}
\fi


\maketitle

\section{Introduction}
\label{sec:introduction}

Metamorphism is the integral transform \(f(u) \mapsto \metamorph{ f}(x,y,b,r)\) defined by~\cites{Kisil21c,Kisil21b}:
\begin{equation}
  \label{eq:metamorphism-intro}
  \begin{split}
\metamorph{ f}(x,y,b,r)
  &= \sqrt[4]{2r^2}
  \int_{\Space{R}{}} f(u)\,
  \exp\left( -\pi  \myhbar \left(
  ( r^{2}-\rmi b)(u-y)^2
   +2 \rmi  (u-y) x
  \right)\right)
  \,\rmd u
  \,.
\end{split}
\end{equation}
Particular values of \(\metamorph{ f}(x,y,b,r)\) encompass many important integral transforms of \(f(u)\), for example:
\begin{itemize}
\item \(\metamorph{ f}(x,0,0,+0)\)  is the Fourier transform.
\item \(\metamorph{ f}(0,y,0,r)\) is the Gauss--Weierstrass(--Hille) transform~\cite{Zemanian67a}.
\item \(\metamorph{ f}(x,y,0,1)\) is the Fock--Segal--Bargmann (FSB) transform~\amscite{Folland89}*{\S~1.6}.
\item \(\metamorph{ f}(x,y,0,r)\) is the Fourier--Bros--Iagolnitzer (FBI) transform~\amscite{Folland89}*{\S~3.3}.
\item \(\metamorph{ f}(x,y,b,1)\)  was used in~\cites{AlmalkiKisil18a,AlmalkiKisil19a} to treat the Schr\"odinger equation.
\item \(\metamorph{ f}(0,0,b,r)\) is  a sort of wavelet transform for
the affine group~\amscite{AliAntGaz14a}*{Ch.~12}.
\item a variant of quadratic Fourier (or Fresnel--Fourier~\cite{Osipov92a}, or integral Gauss~\cite{Neretin11a}, or linear canonical~\cite{HealyKutayOzaktasSheridan16}, etc.) transform.
\end{itemize}
However, the metamorphism is more than just a formal recombination of classical transforms.  For example, it can be used as a sort of transmutation~\cite{KravchenkoSitnik20a,ShishkinaSitnik20a} to reduce the order of  a differential equations~\cite{Kisil21c}: e.g. a second order differential equation can be transformed to a first order admitting a straightforward solution and transparent geometrical structure~\cites{AlmalkiKisil18a,AlmalkiKisil19a}.

Basic properties of the metamorphism can be verified by direct calculations---the path which was intentionally chosen to reduce the amount of prerequisites in the introductory paper~\cite{Kisil21c}.  
Yet, a genuine origin of metamorphism is a covariant transform related to the Schr\"o\-din\-ger--Jacobi group~\cites{Folland89,Berndt07a} as was already presented in the Jupyter notebooks~\cite{Kisil21b} with respective symbolic computations.

This paper systematically utilises the group theory and covariant transform technique to reinstall the metamorphism transform from a scratch. Furthermore, some sister integral transforms are appearing as well. The paper can be seen as a readable narrative  to a Jupyter notebook~\cite{Kisil21b}, which will be frequently referred here to replace some boring calculations. Our main result is a characterisation of the metamorphism  image space in Thm.~\ref{th:metamorphism-image-characterisation}.

We made this paper as accessible as possible. Its reading does not require an advanced knowledge of group representations and the theory of covariant transform. We provide most of required information with further references to more detailed presentations if needed.

In Sect.~\ref{se:groups} we introduce several groups: the Heisenberg, $\SL$, affine, Schr\"o\-din\-ger, and finally our main object---the group SSR. Essential relations between those groups are presented as well. We describe some (not all) induced representations of the group SSR in Sect.~\ref{se:induced-representations}.  The corresponding covariant transform and its properties are described in Sect.~\ref{se:covariant-transform}. Finally, we  connect a selection of a fiducial vector with the properties of the image space of covariant transform in Sect.~\ref{sec:image-space-covar-trans}. In particular, the metamorphism is defined as the covariant transform with a remarkable fiducial vector---the Gaussian. Covariant transforms with some other mentioned fiducial vectors are still awaiting their investigation.

\section{Heisenberg, $\SL$, affine, Schr\"o\-din\-ger and SSR groups}
\label{se:groups}

We start from a brief account of groups involved in the consideration.  An element of the  one-dimensional Heisenberg group \(\mathbb{H}\)~\cites{Folland89,Kisil10a,Kisil19b} will be denoted by \(( s, x, y)\in \Space{R}{3}\). The group law on  \(\mathbb{H}\) is defined as follows: 
\begin{displaymath}
  (s, x, y)\cdot (s ',\, x ',\, y' ) =( s+ s'+ \half{\omega}( x, y; x', y'),\, x+\, x', \,y+ \, y'),
\end{displaymath}
where
\begin{equation}
  \label{eq:symplectic-form}
  {\omega} ( x, y; x', y')= xy'- x'y
\end{equation}
is the symplectic form~\amscite{Arnold91}*{\textsection{}41} on \(\Space{R}{2}\).
The identity element in \(\Space{H}{}\) is \(( 0, 0, 0)\), and the inverse of \(( s, x, y)\) is \(( -s, -x, -y)\).

There is an alternative form of \(\Space{H}{}\) called the polarised  Heisenberg group \(\Space[p]{H}{}\) with the group law~\citelist{\amscite{Folland89}*{\textsection{}1.2} \cite{AlameerKisil21a}}
\begin{displaymath}
  (s,\, x, \,y)\cdot (s ',\, x ',\, y' )= (\, s+ s' + x y',\, x+ x ',\, y+ y ').
\end{displaymath}
and the group isomorphism \(\theta : \Space{H}{}\rightarrow \Space[p]{H}{}\) given by
\begin{displaymath}
  \Theta: ( s, x, y)\rightarrow ( s + \half xy,\, x,\, y).
\end{displaymath}

The special linear group \(\SL\) is the group of \(2\times 2\) matrices with real entries and the unit determinant~\cites{Lang85,Kisil12a}. The group law on  \(\SL\) coincides with the matrix multiplication.
A matrix \(A\in\SL\) acts on vectors in \(\Space{R}{2}\) by a symplectomorphism, i.e. an 
automorphisms of the symplectic form \(\omega\)~\eqref{eq:symplectic-form}:
\begin{displaymath}
  \omega(A(x,y);A(x',y'))=\omega(x,y;x',y').
\end{displaymath}
Therefore, the transformation \(\theta_{A}: \Space{H}{}\rightarrow \Space{H}{}\) 
\begin{displaymath}
  \theta_{A}: ( s, x, y)\rightarrow ( s, A( x, y))
\end{displaymath}
is an automorphism of \(\Space{H}{}\)~\amscite{Folland89}*{\textsection{}1.2}. 
 The corresponding polarised automorphism \(\theta_{A}^{p} = \Theta \circ \theta_A \circ \Theta^{-1}: \Space[p]{H}{}\rightarrow \Space[p]{H}{}\) is
 \begin{displaymath}
   \theta_{A}^{P}( s, x, y)=
  \left( s + \half ( ac x^{2} +2bc xy + bd y^{2}), ax + by, cx+dy\right),
\end{displaymath}
where \(A=\begin{pmatrix} a&b\\c&d\end{pmatrix}\).

Upper-triangular matrices in \(\SL\) with positive diagonal entries form a subgroup  \(\Space{A}{}\). We parameterise it by  pairs \((b, r) \in \Space[+]{R}{2}\) with \(b\in\Space{R}{}\) and  \(r>0\) as follows:
\begin{equation}
  \label{eq:affine-group}
  \begin{pmatrix}1&b\\0&1\end{pmatrix}\,
  \begin{pmatrix}r&0\\0&1\slash r\end{pmatrix}\,
  =\begin{pmatrix}r&b\slash r \\0&1\slash r \end{pmatrix}.
\end{equation}
The subgroup is isomorphic to the affine group of the real line also known as the \(ax+b\) group~\cite{Kisil12d}.


For a group acting by automorphism on another group we can define their semi-direct product. The model case is the affine group itself, where dilations act as automorphisms of shifts. Formally, let \(G\) and \(H\) be two groups and assume \(\theta : H\rightarrow Aut(G)\), where \(\theta_h\) is an automorphism of \(G\) corresponding to \(h\in H\). The semi-direct product of \(G\) by \(H\) denoted by \(G\rtimes H\) is the Cartesian product of \(G\times H\) with the group law
\begin{equation}
  \label{eq:semiderect-law}
  ( g_{1}, h_{1})\,\cdot ( g_{2}, h_{2})= ( g_{1} \theta_{h_{1}}(g_{2}), h_{1} h_{2}),
\end{equation}
where \(( g_{1}, h_{1}),\, ( g_{2}, h_{2})\in G \times H\).

The semidirect product of the Heisenberg group and \(\SL\) is called Schr\"o\-din\-ger group \(\mathbb{S}\), which is the group of symmetries of the Schr\"o\-din\-ger equation~\cites{Niederer72a,KalninsMiller74a} and parabolic equations~\cite{Wolf76a} with applications in optics~\cites{ATorre08a,ATorre10a}. In the context of number theory it is also known as the Jacobi group~\cite{Berndt07a}. 

Our main object here is the group \(\Space{G}{}  \coloneqq  \Space[p]{H}{} \rtimes \Space{A}{}\), which is the semi-direct product  of the Heisenberg group \(\Space[p]{H}{}\) and the affine group \(\Space{A}{}\)~\eqref{eq:affine-group} acting by symplectic automorphism of \(\Space[p]{H}{}\). Thus, \(\Space{G}{}\) is  a subgroup of the Schr\"o\-din\-ger group. It can be also called shear-squeeze-rotation (SSR) group~\cite{Kisil21b} by three types of transformations of Gaussian coherent states. A subgroup of \(\Space{G}{}\) without squeeze (i.e. \(r=1\) in~\eqref{eq:affine-group}) is called the shear group and it was used in a similar context in~\cites{AlmalkiKisil18a,AlmalkiKisil19a}. This nilpotent step 3 group is also known as the Engel group~\cite{Chatzakou22a}.

Let \(( s, x, y,  b, r)\in \Space{G}{}\) 
where \( ( s, x, y) \in \Space[p]{H}{}\) and \(( b, r) \in \mathbb{A}\).
Explicitly the group law~\eqref{eq:semiderect-law} on \(\Space{G}{}\) is~\cite{Kisil21b} 
\begin{align*}
  ( s, x, y, b, r) \cdot ( s', x', y', b', r')
  &=( s+ s'+ x  {r}^{-1} y' - \half\, b \,{({r}^{-1} {y'})}^{2},\\
  & x + r x'- b {r}^{-1} y',\, y +  {r}^{-1} y', \, b + b' {r}^{2},\, r r').
\end{align*}

There is a convenient matrix realisation of \(\Space{G}{}\)~\cite{Kisil21b}
\begin{displaymath}
  (  s, x, y, b, r)=
  \begin{pmatrix}
    1&-yr&({x+ b y})/{r}&2s-yx\\
    0&r&-b/r&x\\
    0&0&1/r&y\\
    0&0&0&1
  \end{pmatrix}.
\end{displaymath} 
The corresponding  solvable Lie algebra \(\algebra{g}\) has a basis \(\{ S, X, Y, B, R\}\), with the following non-vanishing commutators:
\begin{equation}
  \label{eq:lie-algebra-commutators}
  [ X, Y]= S,\quad [ X, R]= -X,\quad [ Y, R]=  Y, \quad
  [ Y, B]= X,\quad  [ R, B]= 2B.
\end{equation}
Clearly, the group \(\Space{G}{}\) is non-commutative.

\section{Induced representations the group $\Space{G}{}$}
\label{se:induced-representations}

In this section we construct several induced representations of the group \(\Space{G}{}\), which are required for our study. First, we recall the general scheme of induced representations. For simplicity, only inductions from characters of subgroups are considered and it is sufficient for our present purposes. For further details and applications of induced representations see~\cites{AliAntGaz14a,Folland16a,Mensky76,Mackey70a,Kisil09e}. 

\subsection{Induced representation from a subgroup character}

Let \(G\) be a group and \(H\) be a subgroup of \(G\). The space \(X= G \slash H\) of the left cosets \(gH\) of the subgroup \(H\) is given by the equivalence relation: \( g\sim g'\) if there exists \( h\in H\) such that \(g = g' h\). We define the natural projection \( \mathbf{p}: G \rightarrow X\) such that \( \mathbf{p}(g)= g H\). 

Let us fix a section \(\mathbf{s}: X\rightarrow G\) such that \(\mathbf{p} \circ \mathbf{s}= I\), where \(I\) is the identity map on \(X\). 
An associated map \(\mathbf{r}: G\rightarrow H\) by 
\begin{equation}
  \label{eq:map-r-formula}
  \mathbf{r}( g)= {\mathbf{s}( \mathbf{p}( g))}^{-1}\cdot g.
\end{equation}
provides  the unique decomposition of the form~\amscite{Kirillov76}*{\textsection{}13.2}
\begin{displaymath}
  g = \mathbf{s} ( \mathbf{p}(g)) \cdot \mathbf{r}(g), \qquad \text{ for any } g \in G.
\end{displaymath}
Thus, \(X\) is a left homogeneous space with the \(G\) action as follows: 
\begin {equation}
  \label{eq:G-space-action}
  g^{-1}: x\rightarrow g^{-1} \cdot x = \mathbf{p}\, ( g^{-1} * \mathbf{s}(x)),
\end{equation}
where \( * \) is the multiplication of \(G\) and \(\cdot\) is the action of \(G\) on \(X\) from the left.

Suppose \(\chi : H\rightarrow \Space{T}{}\) be a character of the subgroup \(H\). 
Let \(\FSpace[\chi]{L}{2}( G)\) be a Hilbert space of functions on \(G\)  with a \(G\)-invariant inner product and the \(H\)-covariance property
\cite{Kisil17a},
\begin{equation}
  \label{eq:covariance-property}
  F( g h) = \Bar{\chi} (h) \, F (g),  \qquad \text{ for all } g\in G, \ h \in H.
\end{equation}
The space \(\FSpace[\chi]{L}{2}( G)\) is invariant under the left regular representation by \(G\)-shifts 
\begin{equation}
  \label{eq:left-regular-action}
  \Lambda (g): F(g')\rightarrow F( g^{-1} g'), \quad \text{ where } g , g' \in G.
\end{equation}
The restriction of \( \Lambda\) to the space \(\FSpace[\chi]{L}{2}( G)\) is called the induced representation from the character \(\chi\).

An equivalent form of the induced representation can be constructed as follows~\cites{Kirillov76,Kisil17a}.  We define a lifting \( \oper{L}^{\chi} : \FSpace{L}{2}(X) \rightarrow \FSpace[\chi]{L}{2}(G) \) as the map
\begin{equation}
  \label{eq:lifting}
    [\oper{L}^{\chi} f] ( g)= \overline{\chi}( \mathbf{r}( g)) \,f (\mathbf{p}( g)).
\end{equation}
The pulling \(\mathcal{P}: \FSpace{L}{2} ^{\chi}(G)\rightarrow \FSpace{L}{2}( X)\)  given by
\begin{equation}
  \label{eq:pulling}
  [\mathcal{P} F] (x)= F (\mathbf{s} ( x)).
\end{equation}
Clearly \(\mathcal{P} \circ \oper{L}^{\chi} = I\) on \(\FSpace{L}{2}(X)\).
From~\eqref{eq:lifting},~\eqref{eq:pulling}, the induced representation \( \rho: \FSpace{L}{2}( X) \rightarrow \FSpace{L}{2}( X)\) is defined by the formula:
\begin{displaymath} 
  \rho_{\chi}( g)  = \mathcal{P} \circ \Lambda ( g) \circ \oper{L}^{\chi},
\end{displaymath}
where \(\Lambda ( g)\) is the left regular representation~\eqref{eq:left-regular-action}.
The representation \(\rho_{\chi}\)  explicitly is
\begin{equation}
  \label{eq:induced-rep-homogen}
  [\rho_{\chi}(g)](x)= \bar{\chi}( \mathbf{r}( g^{-1} \,\mathbf{s} (x))) \, f( g^{-1}\cdot x),
\end{equation}
where \(g\in G\) and \(x\in X\) and  \(g^{-1}\cdot x\) is defined by~\eqref{eq:G-space-action}. For a \(G\)-invariant measure \(\mu\) on \(X\) the representation~\eqref{eq:induced-rep-homogen} is unitary on the space \(\FSpace{L}{2}(X,\mu)\)

\subsection{Derived representations}

In this subsection \(G\) is a Lie group with the corresponding Lie algebra \(\algebra{g}\).
Let \(\rho\) be a representation of \(G\) in a Hilbert space  \(\mathcal{H}\), the derived representation of \(X\in \algebra{g}\) denoted as \(\rmd\rho^{X}\) is given by 
\begin{equation}
  \label{eq:derived-representation}
  \ifspringer
\else
  \rmd{\rho}^{X} \phi = \left. \frac{\rmd\ }{\rmd t}\rho(\exp(t X)) \phi \right|_{t=0},
  \fi
\end{equation}
where the vector \(\phi \in \mathcal{H}\) is such that  the vector-function \(g \rightarrow \rho (g) \phi\) is infinitely-differentiable for any \(g \in G\). These vectors are called smooth and constitute a linear subspace, denoted \(\mathcal{D}^{\infty}\), of \(\mathcal{H}\) which is dense in \(\mathcal{H}\). It is easy to show that \(\mathcal{D}^{\infty}\) is invariant under \(\rho (g)\)~\amscite{Lang85}*{\textsection{}6.1}.
If \(\mathcal{H}\) is  \(\FSpace{L}{2}(\Space{R}{n})\) then the space \(D^{\infty}\) contains the Schwartz space, which is a dense subspace of \(\FSpace{L}{2}(\Space{R}{n})\). 

Also, we define the Lie derivative  \(\mathcal{L}^{X}\) for \(X \in \algebra{g}\) as the derived right regular representation~\amscite{Lang85}*{\textsection{}6.1}, that is
\begin{equation}
  \label{eq:lie-derivative}
  \ifspringer
\else
  [\mathcal{L}^{X} F](g)= \left. \frac{d\ }{dt} F ( g\, \exp(t X))\right|_{t=0},
  \fi
\end{equation}
for any differentiable function \(F\) on \(G\).

\subsection{Left regular representation of group $\Space{G}{}$}
 The left and right invariant Haar measures of the group \(\Space{G}{}\) are given by 
\begin{displaymath} \rmd_l( s, x, y, b, r) =  \rmd s\, \rmd  x \, \rmd  y\, \rmd  b \, \frac{\rmd  r}{{r}^{3}},\end{displaymath}
\begin{displaymath} \rmd_r( s, x, y, b, r) =  \rmd s\, \rmd  x \, \rmd  y \, \rmd  b\, \frac{\rmd r}{r}.\end{displaymath}
Thus, the group \(\Space{G}{} \) is non-unimodular with the modular function  \(\Delta (s, x, y, b, r)=\textstyle \frac{1}{r^{2}}\).

We extend the action~\eqref{eq:left-regular-action} of \(\Space{G}{}\) on itself by left shifts to the left regular unitary representation on the linear space of functions \(\FSpace{L}{2}(\Space{G}{},d_l)\):
\begin{equation}
  \label{G-left-regular}
  \begin{split}
\lefteqn{[ \Lambda ( s, x, y, b, r) F]( s', x', y', b', r')  =F( s'- s + x (y' - y) - \half b {( y' - y)}^{2},}&\qquad{}\\
&\qquad \qquad \qquad \qquad \qquad \frac{1}{r} (x' -x) + \frac{b}{r} (y'- y),\, r (y'- y), \, \frac{1}{r^{2}} ( b'- b),\, \frac{r'}{r}),
\end{split}
\end{equation}
where \(( s, x, y, b, r)\), \(( s', x', y', b', r') \in \Space{G}{}\).

This representation is reducible, i.e. there are \(\Lambda\)-invariant proper subspaces in \(\FSpace{L}{2}(\Space{G}{},d_l)\). In particular, many types of induced representations of \(\Space{G}{}\) are realised as restrictions of the left regular representations~\eqref{G-left-regular} to some subspaces with a covariance property~\eqref{eq:covariance-property}.
We describe here two of them---called the quasi-regular type representation and the Schr\"o\-din\-ger type representation---together with equivalent forms on the respective homogeneous spaces.

\subsection{Quasi-regular representation of the group $\Space{G}{}$}
Let
\begin{displaymath}
  Z=\{( s, 0, 0, 0, 1), s\in \Space{R}{}\}
\end{displaymath}
be the centre of the group \(\Space{G}{}\). The space of left cosets \(X= \Space{G}{}\slash Z\) can be parametrised by
\begin{displaymath}
  \Space[+]{R}{4} = \{ (x, y, b, r) \in \Space{R}{4}: \  r > 0\}.
\end{displaymath}
Consider the natural projection and the section maps
\begin{align}
  \nonumber 
  \mathbf{p}( s, x, y, b, r)&\rightarrow ( x, y, b, r),\\
  \label{eq:map-s-center}
  \mathbf{s}( x, y, b, r) &\rightarrow( 0, x, y, b, r).
\end{align}
We calculate the respective map \( \mathbf{r}\)~\eqref{eq:map-r-formula} as follows
\begin{displaymath}
  \begin{split}
    \mathbf{r} ( s, x, y, b, r)&=  \mathbf{s}( \mathbf{p}( s, x, y, b, r))^{-1} ( s, x, y, b, r)\\
    &=( s, 0, 0, 0, 1).
  \end{split}
\end{displaymath}
Let \(\chi_{\hbar} : Z\rightarrow \Space{T}{}\) be an unitary character of \(Z\):
\begin{displaymath} \chi_{\hbar} ( s, 0, 0, 0, 1) = \rme^{ 2 \pi \rmi  \hbar s},\end{displaymath}
defined by a parameter \( \hbar \in \Space{R}{}\).  In quantum mechanical framework \(\hbar\) is naturally associated to the Planck constant~\cites{Folland89,Kisil02e,Kisil09e,Kisil17a}. 
The corresponding induced representation \(\tilde{\rho} : \FSpace{L}{2}(\Space[+]{R}{4})\rightarrow \FSpace{L}{2}(\Space[+]{R}{4})\) is~\cite{Kisil21b} 
\begin{equation}
  \label{eq:quasi-regular}
  \begin{split}
    \lefteqn{[\tilde{\rho} ( s, x, y, b, r)f]( x', y', b', r')=
      \rme^{ 2 \pi \rmi \hbar( s + x ( y' - y) -  b {( y' -  y)}^{2}/2)}}&\qquad\\
    & \qquad  \qquad  \qquad  \qquad {} \times f( \textstyle\frac{1}{r}( x'-  x) +\textstyle \frac{b}{r} ( y'- y), r ( y'- y), \textstyle\frac{1}{r^{2}}(b' - b ),  \textstyle\frac{r'}{r}).
  \end{split}
\end{equation}
It is called the quasi-regular type representation on \(\FSpace{L}{2}(\Space[+]{R}{4})\). One can check that \(\tilde{\rho}\) is unitary and we will discuss its reducibility below.

\subsection{Schr\"o\-din\-ger type representation of the group $\Space{G}{}$}

Let
\begin{displaymath}
H_{1}=\{( s, x, 0, b, r),\, s,\, x, b \in \Space{R}{}, r \in \mathbb{R_{+}}\}  
\end{displaymath}
be a subgroup of \(\Space{G}{}\), which is a semidirect product of a maximal abelian subgroup of \(\Space{H}{}\) and the affine group \(\Space{A}{}\). The space of the left cosets \( \Space{G}{}\slash H_{1}\) is parameterized by \(\Space{R}{}\). 
We define the natural projection  \(\mathbf{p}: \Space{G}{} \rightarrow \Space{R}{}\) and a section map \(\mathbf{s}: \Space{R}{} \rightarrow \Space{G}{} \) by
\begin{align*}
  \mathbf{p}( s, x, y, b, r) &= y,  \\
  \mathbf{s}( y) & = ( 0, 0, y, 0, 1).
\end{align*}
The respective map \( \mathbf{r}\)~\eqref{eq:map-r-formula} is
\begin{displaymath}
  \begin{split}
    \mathbf{r}( s, x, y, b, r)&=\mathbf{s}( \mathbf{p}( s, x, y, b, r))^{-1} ( s, x, y, b, r)\\
    &=( s, x, 0, b, r) .
  \end{split}
\end{displaymath}  
Let \({\chi}_{\hbar \lambda} : H_{1} \rightarrow \Space{T}{}\) be a character \(H_{1}\)
\begin{displaymath}
  \chi_{ \hbar \lambda} ( s, x, 0, b, r) =   \rme^{ 2\pi \rmi \hbar s} \, r^{ \lambda +\frac{1}{2}},
\end{displaymath}
where \(\hbar\in \Space{R}{}\), \(\lambda \in i{\Space{R}{}}\). For simplicity, we will consider here the case of \(\lambda =0\) only. The induced representation on \( \FSpace{L}{2}(\Space{R}{})\) is~\cite{Kisil21b}
\begin{equation}
  \label{eq:schrodinger-type}
  [\rho ( s, x, y, b, r)f](  u)=\sqrt{r} \, \rme^{ 2\pi \rmi \hbar ( s + x( u - y) -   b {( u - y)}^{2}/2)}  \,  f(  r \,(u - y)).
\end{equation}

\begin{rem}
The structure of this representation can be illuminated through its restrictions to the following subgroups:
\begin{itemize}
\item The affine group \(\mathbb{A}\), i.e. the substitution \(s=x=y=0\). The restriction is the co-adjoint representation of the affine group~\citelist{\amscite{Folland16a}*{\textsection{}6.7.1} \cite{Kisil12d}}:
  \begin{displaymath}
    [\rho ( 0, 0, 0, b, r)f](  u)= \sqrt{r}\, \rme^{ \pi \rmi \hbar \,  b\, { u}^{2} } \,   f(  r\, u).
  \end{displaymath}
  Through the Fourier transform it is unitary equivalent to the quasi-regular representation of the affine group, which is the keystone of the wavelet theory and  numerous results in complex and harmonic analysis~\cite{Kisil12d}.

\item  The Heisenberg group, that is \(r=1\) and \(b=0\).  The restriction is the celebrated Schr\"o\-din\-ger representation~\cites{Folland89,Kisil17a}:
  \begin{displaymath}
    [\rho ( s, x, y, 0, 1)f](  u)=\rme^{ 2\pi \rmi \hbar ( s + x( u - y)) }  \,  f(u - y),
  \end{displaymath}
  which plays the crucial r\^ole in quantum theory.  
    
\item The third subgroup is the Gabor group with  \(b=0\). The representation is
  \begin{displaymath}
    [\rho ( s, x, y, r, 0)f]( u)= \rme^{ 2\pi \rmi \hbar ( s + x( u - y))} r^{ \frac{1}{2} }\,  f(  r \,(u - y)).
  \end{displaymath}
 It is involved in Gabor analysis and Fourier--Bros--Iagolnitzer (FBI) transform~\amscite{Folland89}*{\textsection{}3.3}.
 
\item Finally, the shear group corresponding to \(r=1\). The restriction is 
  \begin{displaymath}
    [\rho ( s, x, y, 1, b)f](  u)= \rme^{ 2\pi \rmi \hbar ( s + x( u - y) -   b {( u - y)}^{2}/2)} f(u - y),
  \end{displaymath}
 It was employed in~\cites{AlmalkiKisil18a,AlmalkiKisil19a} to reduce certain quantum Hamiltonians to first-order differential operators.
\end{itemize}
\end{rem}

In view of the mentioned connections, we call representation~\eqref{eq:schrodinger-type} as Schr\"o\-din\-ger type representation. It is irreducible since its restriction to the Heisenberg group coincides with the irreducible Schr\"o\-din\-ger representation~\cites{Folland89, Kisil17a}.  

The derived  representation~\eqref{eq:derived-representation} of the Schr\"o\-din\-ger type representation~\eqref{eq:schrodinger-type} is
\begin{align}
  \nonumber 
  \rmd{\rho}^{X}&= 2 \pi \rmi \hbar u I, &
                                        \rmd{\rho}^{B}&= - \pi \rmi  \hbar u^{2} I,\\
  \label{eq:schrodinger-derived}
  \rmd{\rho}^{Y}&=  -\frac{\rmd\ }{\rmd u},&
                                    \rmd{\rho}^{R}&= \half  I+u\, \frac{\rmd\ }{\rmd u},\\
  \nonumber 
 \rmd{\rho}^{S}&= 2 \pi \rmi \hbar I.
\end{align}
It is easy to check that the above sets of operators 
\eqref{eq:schrodinger-derived} represents commutators~\eqref{eq:lie-algebra-commutators} of the Lie algebra \(\algebra{g}\) of the group \(\Space{G}{}\).

\section{Covariant transform}
\label{se:covariant-transform}

The covariant transform plays a significant r\^ole in various fields of mathematics and its applications~\cites{Perelomov86, Berezin86, AliAntGaz14a, Folland89, Kisil11c, Kisil17a, Kisil12d}. We present here some fundamental properties of the covariant transform which have implications for the metamorphism transform.

\subsection{Induced covariant transform}

Let \(G\) be a group and let \(\rho\) be a unitary irreducible representation of the group \(G\) in a Hilbert space \(\mathcal{H}\).  For a fixed unit vector \(\phi \in \mathcal{H}\), called here a fiducial vector (aka vacuum vector, ground state, mother wavelet, etc.), the covariant transform \(\oper{W}_{\phi}: \mathcal{H}\rightarrow  \FSpace{L}{}(G)\)
\ifspringer
\else
is~\citelist{\cite{AliAntGaz14a}*{\textsection{}8.1} \cite{Berezin86} \cite{Perelomov86}}
\fi
\begin{equation}
  \label{eq:covariant-transform}
  [ {\oper{W}_{\phi}}f](g)= \langle f, \rho(g) \phi \rangle, \qquad
  \text{ where } f\in \mathcal{H} \text{ and } g\in G.
\end{equation}
Here \(\FSpace{L}{}(G)\) is a certain linear space of functions on \(G\) usually linked to some additional conditions. The common focus is on \(\FSpace{L}{}(G)=\FSpace{L}{2}(G, \rmd\mu)\)---the square integrable functions with respect to a Haar measure \(\rmd\mu\), for this we need a square-integrable representation \(\rho\) and an admissible fiducial vector \(\phi\)~\cite[Ch.~8]{AliAntGaz14a}. However, many topics in analysis prompt to study other situations, e.g. related to Hardy-type invariant functionals~\cites{Kisil12d,Kisil13a}, Gelfand triples~\cite{FeichtingerGrochenig88a},  Banach spaces~\cite{Kisil98a}, etc.

The main property of~\eqref{eq:covariant-transform} is that \({\oper{W}_{\phi}}\) intertwines the representation \(\rho\) on \(\mathcal{H}\) and the left regular action  \(\Lambda\)~\eqref{eq:left-regular-action}  on \(G\):
\begin{equation}
  \label{eq:covariant-intertwining}
  \oper{W}_{\phi} \circ \rho( g) = \Lambda (g) \circ \oper{W}_{\phi},\quad \text{ for all  } g\in G.
\end{equation}

A representation \(\rho\) is square-integrable if for some \(\phi \in \mathcal{H}\), the map \(\oper{W}_{\phi}: \mathcal{H} \rightarrow \FSpace{L}{2}(G, dg)\) is unitary for a left Haar measure  \(dg\) on \(G\). Some representations are not square-integrable, but can still be treated by the following modification of covariant transform by Perelomov~\cite{Perelomov86}. Let \(H\) be a closed subgroup of the group \(G\) and the corresponding homogeneous space is \(X= G\slash H\).  Let for some character \(\chi\) of \(H\)  a fiducial vector  \(\phi \in \mathcal{H}\) is a joint eigenvector
\begin{equation}
  \label{eq:eigenvector-subgroup}
  \rho (h) \, \phi = \chi (h) \phi, \qquad \text{for all } h \in H.
\end{equation}
Then, the respective covariant transform satisfies the covariant property, cf.~\eqref{eq:covariance-property}:
\begin{displaymath}
    [\oper{W}_{\phi} f](g h)
    = \overline{\chi}(h) [\oper{W}_{\phi} f] (g).
\end{displaymath}
Thus, the image space of \(\oper{W}_{\phi}\) belongs to the induced representation by the character \(\chi\) of the subgroup \(H\). This prompts to adopt the covariant transform to the space of function on the homogeneous space \(X=G/H\).
To this end, let us fix a section \(\mathbf{s}: X \rightarrow G\) and a fiducial vector \(\phi\in \mathcal{H}\) satisfying~\eqref{eq:eigenvector-subgroup}. The induced covariant transform from the Hilbert space \(\mathcal{H}\) to a space of functions \(\FSpace{L}{\phi}(X)\) is
\begin{displaymath}
  [ {\oper{W}_{\phi}}f](x)= \langle f, \rho(s(x)) \phi \rangle,\quad
  \text{ where }  x\in X.
\end{displaymath}
Then, the induced covariant transform intertwines \(\rho\) and \(\tilde{\rho}\)---an induced representation from the character \(\chi\) of the subgroup \(H\), cf.~\eqref{eq:covariant-intertwining}:
\begin{equation}
  \label{eq:covariant-intertwining-induced}
    W^{\rho}_{\phi}\circ \rho(g) = \tilde{\rho}(g) \circ  W^{\rho}_{\phi}, \qquad \text{ for all } g\in G.
 \end{equation}
In particular, the image space \(\FSpace{L}{\phi}( G\slash H)\) of the induced covariant transform is invariant under \(\tilde{\rho}\). Induced covariant transforms for the Heisenberg group~\cite{Kisil09e} and the affine group~\cite{Kisil12d} are the most familiar examples.

\subsection{Induced covariant transform of the group $\Space{G}{}$}
On the same way as above, we can calculate the induced covariant transform of \(\Space{G}{}\).
Consider the subgroup \(Z\) of \(\Space{G}{}\), which is \(Z= \{( s, 0, 0, 0, 1), s\in \Space{R}{}\}\). For the Schr\"o\-din\-ger type representation~\eqref{eq:schrodinger-type}, any function  \(\phi \in \FSpace{L}{2}(\Space{R}{})\) satisfies the eigenvector condition \(\rho( s, 0, 0, 0, 1) \phi= \rme^{2 \pi \rmi \hbar s} \phi\) with the character \(\chi ( s, 0, 0, 0, 1)= \rme^{-2 \pi \rmi \hbar s}\), cf.~\eqref{eq:eigenvector-subgroup}. Thus, the respective homogeneous space is \(\Space{G}{}\slash Z \simeq \Space[+]{R}{4}\) and we take the above section  \(\mathbf{s}: \Space{G}{}\slash Z \rightarrow \Space{G}{}: \ \mathbf{s}( x, y, b, r)= ( 0, x, y, b, r)\)~\eqref{eq:map-s-center}. Then, the induced covariant transform is
\begin{equation}
  \label{eq:covariant-tr-schrodinger-type}
    \begin{split}
      [ \oper{W}_{\phi} f]( x, y, b, r)&= \langle f, \rho(\mathbf{s} ( x, y, b, r)) \phi \rangle \\ 
      &= \langle f,\rho( 0,x, y, b, r) \phi \rangle\\ 
      &= \int_{\Space{R}{}}{f (u)\, \overline{\rho( 0, x, y, b, r) \, \phi ( u)}} \, \rmd  u\\
      &= \int_{\Space{R}{}}{ f ( u) \,\rme^{- 2 \pi \rmi \hbar (x ( u -y) - b (u - y)^{2}/2)} \,r^{ \frac{1}{2}} \, \overline\phi (r( u- y))} \, \rmd  u\\ 
      &=\sqrt{r} \int_{\Space{R}{}}{ f ( u)\, \rme^{-2 \pi \rmi \hbar (x ( u -y) - b (u - y)^{2}/2)}\, \overline\phi ( r( u- y))} \, \rmd  u.
    \end{split}
\end{equation}
From~\eqref{eq:covariant-intertwining-induced}, \(\oper{W}_{\phi}\) intertwines the Schr\"o\-din\-ger type representation~\eqref{eq:schrodinger-type} with quasi-regular~\eqref{eq:quasi-regular}.

The last integral in~\eqref{eq:covariant-tr-schrodinger-type} is a composition of five unitary operators \(\FSpace{L}{2}(\Space{R}{2})\rightarrow \FSpace{L}{2}(\Space{R}{2}) \) applied to a function \(F(y,u)=f(y)\overline {\phi}(u)\) in the space \(\FSpace{L}{2}(\Space{R}{})\otimes  \FSpace{L}{2}(\Space{R}{}) \simeq \FSpace{L}{2}(\Space{R}{2})\):
\begin{enumerate}
 \item The unitary operator \(R: \FSpace{L}{2}(\Space{R}{2}) \rightarrow  \FSpace{L}{2}(\Space{R}{2})\) based on the dilation
   \begin{displaymath}
     R : F( y, u)\rightarrow \sqrt{r}\, F( y, ru), \qquad \text{    where } r > 0.    
   \end{displaymath}

 \item The change of variables  \(T: \FSpace{L}{2}(\Space{R}{2}) \rightarrow  \FSpace{L}{2}(\Space{R}{2})\)
      \begin{displaymath}
        T: F( y, u)\rightarrow F(u, u-y).
      \end{displaymath}
    \item The operator of multiplication by an unimodular function \(\psi_{b}( x, y)= \rme^{\pi \rmi \hbar b (u -y)^{2}}\), 
      \begin{displaymath}
        M_{b}: F( y,u) \rightarrow  \rme^{\pi \rmi \hbar b (u -y)^{2}} \,F ( y, u), \qquad
        \text{ where } b\in \Space{R}{}. 
      \end{displaymath}
    
    \item The partial Fourier transform \(u\rightarrow x\) in the second variable 
      \begin{displaymath}
        [\mathcal{F}_{2} F]( y, x)= \int_{\Space{R}{}}{ F( y, u)\, \rme^{- 2 \pi \rmi \hbar x  u}\, \rmd u}.
      \end{displaymath}
    \item The multiplication \(M\) by the unimodular function \(\rme^{2 \pi \rmi \hbar x y}\).
\end{enumerate}
Thus, we can write \(\oper{W}_{\phi}\) as
\begin{equation}
  \label{eq:covarinat-as-composition}
  [\oper{W}_{\phi}f]( x, y, b, r) =  [\left(M \circ \mathcal{F}_{2} \circ M_{b} \circ T \circ R\right)\,F ]( x, y),
\end{equation}
and obtain
\begin{prop}
  \label{pr:unitarity-G-covariant}
  For a fixed \(r_{0}\in\Space[+]{R}{}\) and \(b_{0}\in \Space{R}{}\), the map \(f \otimes \overline{\phi}\rightarrow [\oper{W}_{\phi} f](\cdot, \cdot, b_0, r_0)\) is a unitary operator from  \(\FSpace{L}{2}(\Space{R}{})\otimes \FSpace{L}{2}(\Space{R}{})\) onto \( \FSpace{L}{2}(\Space{R}{2})\).
\end{prop}

Also, the induced covariant transform preserves the Schwartz space, that is, if \(f, \phi \in \mathcal{S}(\Space{R}{})\) then \(\oper{W}_{\phi}f ( \cdot, \cdot, b_0, r_0)\in \mathcal{S}(\Space{R}{2})\). This is because the \(\mathcal{S}(\Space{R}{2})\) is invariant under the all five above components of \(\oper{W}_{\phi}f\) in~\eqref{eq:covarinat-as-composition}. 

Note, that the induced covariant transform~\eqref{eq:covariant-tr-schrodinger-type} does not define a square-integrable function on \(\Space{G}{}\slash Z\sim \Space[+]{R}{4}\). To discuss unitarity we need to introduce a suitable inner product. In general we can start from a probability measure \(\mu\) on \(\Space[+]{R}{2}\), that is \(\int_{\Space[+]{R}{2}} \,\rmd \mu =1\). Then we define the inner product 
\begin{equation}
  \label{scalar-product-mu}
  \langle f, g \rangle_{\mu} = \int_{\Space[+]{R}{4}}{f ( x, y, b, r) \, \overline{ g( x, y, b, r)}\, \frac{\hbar \, \rmd  x \, \rmd y \, \rmd \mu(b,r)}{\sqrt{2 r}},}
\end{equation}
for \(f, g \in \FSpace{L}{\phi}(\Space[+]{R}{4})\). The factor \(\hbar\) in the measure \(\textstyle\frac{\hbar \, \rmd  x \, \rmd y}{\sqrt{2 r_{0}}}\) makes it dimensionless, see discussion of this in~\cites{Kisil02e,AlmalkiKisil18a}.
Important particular cases of probability measures parametrised by \((b_0, r_0)\in \Space[+]{R}{2}\) are
\begin{equation}
  \label{eq:dirac-measure}
  \rmd \mu_{(b_0,r_0)}(b,r) = \delta(b-b_0)\, \delta(r-r_0) \,\rmd b \,\rmd r \, ,
\end{equation}
where \(\delta(t)\) is the Dirac delta. The respective inner products becomes:
\begin{equation}
  \label{scalar-product-b0-r0}
  \langle f, g \rangle_{(b_0, r_0)} = \int_{\Space{R}{2}}f ( x, y, b_0, r_0) \, \overline{ g( x, y, b_0, r_0)}\, \frac{\hbar \, \rmd  x \, \rmd y}{\sqrt{2 r_{0}}}\,.
\end{equation}
From now on we consider \(\FSpace{L}{\phi}(\Space[+]{R}{4})\) as a Hilbert space with the inner product~\eqref{scalar-product-mu}  or specifically \eqref{scalar-product-b0-r0}. The respective norms are denoted by \(\| \cdot \|_{\mu}\) and \(\| \cdot \|_{(b_0, r_0)}\).

Using the above inner product, we can derive from Proposition~\ref{pr:unitarity-G-covariant} the following orthogonality relation:
\begin{cor}
  Let \(f, g, \phi, \psi \in \FSpace{L}{2}(\Space{R}{})\), then  
  \begin{displaymath}
  \langle \oper{W}_{\phi} f, \oper{W}_{\psi} g\rangle_{\mu}= \langle f, g \rangle  \, \overline{\langle \phi, \psi \rangle },
\end{displaymath}
for any probability measure \(\mu\), in particular~\eqref{eq:dirac-measure} with fixed  \((b_0, r_0)\in \Space[+]{R}{2}\). 
 \end{cor}
\begin{cor}
  Let \(\phi \in \FSpace{L}{2}( \Space{R}{})\) have a unite norm. Then, the induced covariant transform \(\oper{W}_{\phi}\) is an isometry from \(\FSpace{L}{2}(\Space{R}{})\) to
  \(\FSpace{L}{\phi}( \Space[+]{R}{4})\) and its inverse is given by the adjoint operator---contravariant transform:
  \begin{equation}
    \label{eq:contravariant-transform-mu}
    f(u) = \int_{\Space[+]{R}{4}}F ( x, y, b, r) \, [\rho( \mathbf{s}(x, y, b, r)) \phi](u)\, \frac{\hbar \, \rmd  x \, \rmd y \, \rmd \mu(b,r)}{\sqrt{2 r}},
  \end{equation}
  where  \(F \in \FSpace{L}{\phi}( \Space[+]{R}{4})\). In particular:
  \begin{equation}
    \label{eq:contravariant-transform}
    f(u) = \int_{\Space{R}{2}} F ( x, y, b_0, r_0) \, [\rho( \mathbf{s}(x, y, b_0, r_0)) \phi](u) \, \frac{\hbar \, \rmd  x \, \rmd y}{\sqrt{2 r_{0}}}.
  \end{equation}
\end{cor}
\begin{proof}
  For \(f\in \FSpace{L}{2}(\Space{R}{})\), we have
  \begin{displaymath}
    \| f\|_{\FSpace{L}{2}(\Space{R}{})} = \| f \otimes \overline{\phi} \|_{\FSpace{L}{2}(\Space{R}{2})}= \| \oper{W}_{\phi}f\|_{\mu},
  \end{displaymath}
  which follows from isometry \( \FSpace{L}{2}(\Space{R}{2}) \rightarrow \FSpace{L}{2}(\Space{R}{2})\) in Prop.~\ref{pr:unitarity-G-covariant}. Then verification of formulae~\eqref{eq:contravariant-transform-mu}--\eqref{eq:contravariant-transform} is a technical exercise.
\end{proof}
A reader may note that~\eqref{eq:contravariant-transform} with \(\phi(u) = \sqrt[4]{2} \rme^{-\pi \hbar\, u^{2}} \) is essentially the inverse Fock--Segal--Bargmann transform. 

\section{Image spaces of the covariant transforms}
\label{sec:image-space-covar-trans}

Clearly, not every function on \(\Space[+]{R}{4}\) is a covariant transform~\eqref{eq:covariant-tr-schrodinger-type} of a function from \(\FSpace{L}{2}(\Space{R}{})\). In this section we discuss the image space of the covariant transform.

\subsection{Right shifts and covariant transform}
Let \(R(g)\) be the right regular representation of the group \(G\), which acts on the functions defined in the group \(G\) as follows:
\begin{displaymath}
  R(g): f( g')\rightarrow f( g' \, g), \qquad \text{ where } g\in G.
\end{displaymath}
In contrast to the intertwining property of the covariant transform for the left regular representation~\eqref{eq:covariant-intertwining}, the right shift satisfies the relation
\begin{align*}
  R(g) [\oper{W}_{\phi}] (g') & =[\oper{W}_{\phi}]( g' \,g) \\
                              & = \langle f, \rho( g'\,g) \phi\rangle  \\
                              & =\langle f ,\rho(g') \rho(g) \phi \rangle \\
                              & =[\oper{W}_{\rho (g)\,\phi} f](g').
\end{align*}
That is, the covariant transform intertwines the right shift with the action of \(\rho\) on the fiducial vector \(\phi\).  Therefore, we obtain the following result, which plays an important r\^ole in exploring the nature of the image space of the covariant transform.

\begin{cor}
  \textup{\cite{Kisil10c}}
  \label{co:analyticity-fiducial}
  Let \(G\) be a Lie group with a Lie algebra \(\algebra{g}\) and \(\rho\) be a representation of \(G\) in  a Hilbert space \(\mathcal{H}\). Let a fiducial vector \(\phi\) be a null-solution, \(A\phi=0\), for the operator \(A= \sum_{j} a_{j}{d\rho}^{X_{j}}\), where \({d\rho}^{X_{j}}\) are the derived representation of some \(X_{j}\in \algebra{g}\) and \(a_{j}\) are constants. Then, for any \(f\in \mathcal{H}\)  the wavelet transform \([W_{\phi}f](g)= \scalar{f}{\rho(g)\phi}\) satisfies 
  \begin{displaymath}
    D(W_{\phi}f)=0, \quad  \text{ where } \quad D= \sum_{j} \overline{a}_{j} \mathcal{L}^{X_{j}}.
  \end{displaymath} 
  Here \(\mathcal{L}^{X_{j}}\) are the left invariant fields (Lie derivatives)~\eqref{eq:lie-derivative} on \(G\) corresponding to \(X_{j}\).
\end{cor}
Illustrative examples are the classical spaces of analytical  functions: the Fock--Segal--Bargmann space and the Hardy space, see~\cites{Kisil11c,Kisil13c} for details. 
\begin{rem}
  \label{rem:analyticity-polynom}
  It is straightforward to extend the result of Cor.~\ref{co:analyticity-fiducial} from a linear combination of elements in the Lie algebra \(\algebra{g}\) to an arbitrary polynomial from the enveloping algebra of \(\algebra{g}\) or even to more general functions/distributions, cf.~\amscite{Kisil11c}*{Cor.~5.8}.
\end{rem}

\subsection{Characterisation of the image space for the group $\Space{G}{}$}

The above Cor.~\ref{co:analyticity-fiducial} can be used to construct covariant transforms with desired properties through purposely selected fiducial vectors. We are  illustrating this for the group \(\Space{G}{}\). First,  we need to compute the Lie derivatives~\eqref{eq:lie-derivative} reduced to the representation space of the quasi-regular representation~\eqref{eq:quasi-regular}, see~\cite{Kisil21b}:
\begin{align} 
\nonumber 
    \mathcal{L}^{X}&= r \partial_{x},  &
                                         \mathcal{L}^{B}&= r^{2}\,\partial_{b},\\
  \label{eq:lie-derivatives}
  \mathcal{L}^{Y} &=\textstyle\frac{1}{r}( -2\pi \rmi \hbar x I- b\,\partial_{x} + \partial_{y}), &
                                                                                                    \mathcal{L}^{R} &= r\,  \partial_{r} ,\\
  \nonumber 
  \mathcal{L}^{S}&= -2 \pi \rmi \hbar I. 
\end{align}
One can check that those Lie derivatives make a representation of the Lie algebra of the group \(\Space{G}{}\)~\cite{Kisil21b}.

Now we are looking for a covariant transform \(\oper{W}_\phi: \FSpace{L}{2}(\Space{R}{}) \rightarrow \FSpace{L}{2}(\Space[+]{R}{4})\) with the image space annihilated by a generic linear combination of Lie derivatives~\eqref{eq:lie-derivatives}. To this end the fiducial vector \(\phi\) shall be a null solution of the following differential operator composed from the derived  Schr\"o\-din\-ger type representation~\eqref{eq:schrodinger-derived}
\begin{equation}
  \label{eq:generic-derivative}
  \begin{split}
\lefteqn{    {\rmd\rho}^{\rmi E_s S + E_{x}X+ \rmi E_y  Y+ \rmi  E_{b} B+ E_{r} R}}\qquad&\\
    &= (E_ru-\rmi E_y) \frac{\rmd\ }{\rmd u} +(\pi  \hbar(E_b  u^{2}+{2 \rmi }  E_x  u - 2   E_s)+\half E_{r} )I.
  \end{split}
\end{equation}
where \(E_s\), \(E_x\), \(E_y\), \(E_b\) and \(E_r\) are arbitrary real coefficients.
This equation has two different solutions depending on a value of \(E_r\). If \(E_r=0\) (which requires \(E_y\neq 0\) for non-trivial operator~\eqref{eq:generic-derivative}) then a generic solution of~\eqref{eq:generic-derivative} is~\cite{Kisil21b}
\providecommand{\genphizero}{\phi_0}
\begin{equation}
  \label{eq:fiducial-Er=0}
    \genphizero (u)= C\,\exp\!\left(\pi \hbar\left(  2\rmi \frac{E_s}{E_y} u+ \frac{E_x}{E_y} u^2 - \rmi \frac{E_b}{3E_y} u^3  \right)\right)\,,
\end{equation}
where  \(E_x< 0\) for square integrability of \(\genphizero\) and the constant \(C\) is determined from the normalisation condition \(\|\genphizero\|_2 = 1\).  We have here a sort of Airy beam~\cite{BerryBalazs79a}, which was employed in~\cite{AlmalkiKisil19a} in the context of the share group, i.e. the absence of \({d\rho}^R\) in~\eqref{eq:generic-derivative}.

For \(E_r\neq 0\) we find the generic solution in the form~\cite{Kisil21b}:
\providecommand{\genphi}{\phi_1}
\begin{equation}
  \label{eq:fiducial-Er-not-0}
  \begin{split}
    \genphi (u)&= C \, {(E_{r} u - \rmi E_y)^{- \frac{ 1}{2}+2 
        { \pi \hbar E_s}/{E_r}  - 
        { \pi \hbar E_y( 2E_{x} E_{r}+ E_{b} E_y)}/{E_{r}^{3}}}}\\
    &\qquad \times {\exp\left(-\pi \hbar\left(\frac{\rmi(  2E_{x}E_{r}+ E_{b} E_y)}  {E_{r}^{2}} u +\frac{ E_{b}}{2 E_{r}} {u}^{2}\right)\right)}
    .
\end{split}
\end{equation}
Again, for \(\genphi\in \FSpace{L}{2}(\Space{R}{})\) we need \(\frac{\hbar E_{b}}{E_{r}}>0\) and a proper normalising constant \(C\).

A detailed study of all arising covariant transforms is still awaiting further work. Here we concentrate on some special aspects which appear in this extended group setting for the most traditional fiducial vector---the Gaussian.  First, we note that it steams from both solutions~\eqref{eq:fiducial-Er=0} and~\eqref{eq:fiducial-Er-not-0}:
\begin{itemize}
\item For \(E_r=0\) letting \(E_s=E_b=0\), \(E_x=-1\) and \(E_y=1\) with \(C=\sqrt[4]{2}\) in \(\genphizero\)~\eqref{eq:fiducial-Er=0} produces 
  \begin{equation}
    \label{eq:gaussin-first}
    \phi(u) = \sqrt[4]{2} \rme^{-\pi \hbar\, u^{2}} \quad \text{ with the identity } \quad   {\rmd \rho}^{ -X+ \rmi  Y}\, \phi = 0,
  \end{equation}
  i.e. \(\phi\) is annihilated by the Heisenberg group part of \(\Space{G}{}\).
\item For \(E_r=1\) substitution of \(E_s= \frac{1}{4 \pi \hbar }\), \(E_{x}=E_y=0\) and \(E_{b}=2\) with \(C=\sqrt[4]{2}\) into the vacuum \(\genphi\)~\(\genphi\)~\eqref{eq:fiducial-Er-not-0} again produces
  \begin{equation}
    \label{eq:gaussian-second}
  \phi(  u)= \sqrt[4]{2} \rme^{-\pi \hbar\, u^{2}} \quad \text{ with the identity } \quad  {\rmd \rho}^{\rmi/(4\pi \hbar) S+ 2\rmi B+ R}\, \phi = 0,
\end{equation}
i.e. \(\phi\) is also annihilated by the affine group part of \(\Space{G}{}\).

\end{itemize}

Let us introduce the covariant transform \(\oper{W}_\phi: \FSpace{L}{2}(\Space{R}{}) \rightarrow \FSpace{L}{2}(\Space[+]{R}{4})\)~\eqref{eq:covariant-tr-schrodinger-type} with the fiducial vector \(\phi\)~\eqref{eq:gaussin-first}--\eqref{eq:gaussian-second}:
\begin{equation}
  \label{eq:metamorphism}
  \begin{split}
    \lefteqn{
      [\oper{W}_\phi
    f](x,y,b,r)
     =\sqrt{r} \int_{\Space{R}{}}{ f ( u)\, \rme^{-2 \pi \rmi \hbar (x ( u -y) - b (u - y)^{2}/2)}\, \overline\phi ( r( u- y))} \, \rmd  u} \qquad &\\
     &=\sqrt[4]{2r^2} \int_{\Space{R}{}} f ( u)\, \rme^{-2 \pi \rmi \hbar (x ( u -y) - b (u - y)^{2}/2)}\, \rme^{ -\pi \hbar r^2( u- y)^2} \, \rmd  u\\
  &= \sqrt[4]{2r^2}
  \int_{\Space{R}{}} f(u)\,
  \exp\left( -\pi  \myhbar \left(
  ( r^{2}-\rmi b)(u-y)^2
   +2 \rmi  (u-y) x
  \right)\right)
  \,\rmd u
  \,.
\end{split}
\end{equation}
That is we obtained a representation of the metamorphism~\eqref{eq:metamorphism-intro} as a covariant transform: \(    \metamorph{ f}  = \oper{W}_\phi f\). Now the notation \(\metamorph{ f}\) can be explained as the double covariant transform for the Heisenberg and the affine groups simultaneously.
The image space \(\FSpace{L}{\phi}(\Space[+]{R}{4})\) of the metamorphism is a subspace of square-integrable functions on \(\FSpace{L}{2}(\Space[+]{R}{4}, \|\cdot\|_{\mu})\), see~\eqref{scalar-product-mu}.

\begin{rem}
  Another feature of the Gaussian as a fiducial vector is that an extension of the group \(\Space{G}{}\) to the full Schr\"o\-din\-ger group does not add a value. Indeed,  the Iwasawa decomposition \(\SL=ANK\)~\citelist{\amscite{Lang85}*{\textsection{}III.1} \amscite{Kisil12a}*{\textsection{}1.1}} represents \(\SL\) as the product of the affine subgroup \(AN\) and the compact subgroup \(K\). Yet, the Gaussian is invariant under the action  of the phase-space rotations produced by \(K\). Thus, we get the same set of coherent states from the actions of the group \(\Space{G}{}\) and the Schr\"o\-din\-ger group.
\end{rem}

From the annihilation property~\eqref{eq:gaussin-first} by the derived representation \({d\rho}^{ -X+ \rmi  Y}\)  and Cor.~\ref{co:analyticity-fiducial} we conclude that \(\mathcal{L}^{-X - \rmi Y}\,   \metamorph{ f} =0\) for any \(f\).  Using~\eqref{eq:lie-derivatives} we find~\cite{Kisil21b}:
\begin{equation}
  \label{eq:Cauchy-first}
  \begin{split}
    \oper{C}_1 &= -\mathcal{L}^{X} - \rmi\mathcal{L}^{Y}\\
    &=   \frac{1}{r}\left({(r^{2} - \rmi b)}\, \partial_x+ \rmi \, \partial_y+2  x \myhbar \pi\, I\right).   
  \end{split}
\end{equation}
The operator \(\oper{C}_1\) is called the first Cauchy--Riemann type operator. Similarly from~\eqref{eq:gaussin-first} we conclude that \(\oper{C}_2   \metamorph{ f} =0\) for the second Cauchy--Riemann type operator~\cite{Kisil21b}:
\begin{equation}
  \begin{split}
    \label{eq:Cauchy-second}
    \oper{C}_2  & = -\frac{\rmi}{4 \pi \hbar }\mathcal{L}^{S} -2 \rmi \mathcal{L}^{B}  + \mathcal{L}^{R} \\
    &=  2  r^{2}\, \partial_b +\rmi  r\, \partial_r - \half \rmi\, I\,. 
  \end{split}
\end{equation}
It is convenient to view operators \(\oper{C}_1\) and \(\oper{C}_2\) as the Cauchy--Riemann operators for the following complexified variables:
\begin{equation}
  \label{eq:complex-variables}
  w=b + \rmi r^2 \quad \text{ and } \quad  z=x+(b + \rmi r^2)y = x+wy \,.
\end{equation}
\begin{rem}
  As was pointed out in~\cite{Kisil13c}, the analyticity conditions \(\oper{C}_1\metamorph{ f} =0\)~\eqref{eq:Cauchy-first} and \(\oper{C}_2\metamorph{ f} =0\)~\eqref{eq:Cauchy-second} are consequences of minimal uncertainty properties of the fiducial vector. The first condition~\eqref{eq:Cauchy-first} follows from the celebrated Heisenberg--Kennard uncertainty relation~\cites{Folland89,Kisil13c}
  \begin{displaymath}
    \Delta_\phi(M) \cdot \Delta_\phi(D) \geq \frac{\myh}{2}
  \end{displaymath}
  for the coordinate \(M=\rmd{\rho}^{\rmi X}\) and momentum \(D=\rmd{\rho}^{\rmi Y}\) observables in the Schr\"o\-din\-ger representation~\eqref{eq:schrodinger-derived}. The second condition~\eqref{eq:Cauchy-second} is due to the similar minimal joint uncertainty of the Gaussian state for the Euler operator \({\rmd \rho}^{1/(4\pi \hbar) S-\rmi R}=-\rmi u \, \frac{\rmd\ }{\rmd u}\) and the quadratic potential \({\rmd \rho}^{\rmi B}= \pi \hbar u^{2} I\).
\end{rem}

Besides the two operators \(\oper{C}_1\) and  \(\oper{C}_2\) which are based on the special properties~\eqref{eq:gaussin-first}--\eqref{eq:gaussian-second} of the Gaussian we can note a couple of polynomial identities in the Schr\"o\-din\-ger type representation of the Lie algebra \(\algebra{g}\). Indeed, using~\eqref{eq:schrodinger-derived} one can check:
\begin{equation}
  \label{eq:second-order-rep-zero}
  \left(\rmd{\rho}^{X}\right)^2+2\, \rmd {\rho}^{S}\, \rmd {\rho}^{B} =0, \quad \text{ and } \quad
  \rmd{\rho}^{X}\, \rmd {\rho}^{Y}+d{\rho}^{Y}\, \rmd {\rho}^{X} + 2\, \rmd {\rho}^{S}\, \rmd {\rho}^{R} =0.
\end{equation}
These relations express the affine subalgebra generators \(B\) and \(R\) through the Heisenberg ones \(X\) and \(Y\). That is related to so-called quadratic algebra concept~\amscite{Gazeau09a}*{\textsection{}2.2.4}. Because operators in~\eqref{eq:second-order-rep-zero} annihilate any function, including the fiducial vector of the metamorphism, Rem.~\ref{rem:analyticity-polynom} implies that the image space \(\FSpace{L}{\phi}(\Space[+]{R}{4})\) is annihilated by by the second-order  differential operators~\cite{Kisil21b}:
\begin{align}
  \label{eq:strctural-1}
  \oper{S}_1 &=\left(\mathcal{L}^{X}\right)^2+2\,\mathcal{L}^{S}\,\mathcal{L}^{B} = r^{2}( 4 \pi \rmi   \myhbar \, \partial_{b} -   \partial_{xx}^{2}) \,;
\intertext{and}
\label{eq:strctural-2}
  \begin{split}
  \oper{S}_2 &= \mathcal{L}^{X}\,\mathcal{L}^{Y}+\mathcal{L}^{Y}\,\mathcal{L}^{X} + 2\,\mathcal{L}^{S}\,\mathcal{L}^{R}\\
  & = -4  \pi \rmi r \myhbar \, \partial_{r}-  2b \, \partial_{xx}^{2}+2\, \partial_{xy}^{2} - 4  \pi \rmi x \myhbar \, \partial_{x}-2 \pi \rmi   \myhbar\, I\,.
\end{split}
\end{align} 
The identities  \(\oper{S}_1\metamorph{ f}  = 0 \) and  \(\oper{S}_2\metamorph{ f}  = 0 \)  are called the first and second structural conditions. 
Their presence is a  notable difference between covariant transforms on the group \(\Space{G}{}\) and the Heisenberg group.

Of course, the list of annihilators is not exhausting and the above conditions are not independent. If \( \oper{S}_1 F =0\) for a function \(F\) satisfying both the Cauchy--Riemann-type operators~\eqref{eq:Cauchy-first}--\eqref{eq:Cauchy-second}   then the function \(F\) has to be in the kernel of \(\oper{S}_2\)~\eqref{eq:strctural-2} as well, see~\cite{Kisil21b}.

It was shown~\cite{Kisil21c} that a generic solution of two differential operators~\eqref{eq:Cauchy-first}--\eqref{eq:Cauchy-second} is:
\begin{equation}
  \label{eq:first-order-generic-complex}
  [\oper{G}f_2](z,w)   \coloneqq
                       \sqrt{r}\, \rme^{-\pi i \myhbar x^2/w}\, f_2 (z, w)\,,
\end{equation}
where \(f_2\) is a holomorphic function of two complex variables \(z\) and \(w\)~\eqref{eq:complex-variables}. Additionally,  the structural condition \(\oper{S}_1\)~\eqref{eq:strctural-1} applied to \(\oper{G}f_2\)~\eqref{eq:first-order-generic-complex} produces a parabolic equation for \(f_2\):
\begin{equation}
    \label{eq:strctural-generic}
    4\pi i h w\partial_wf_2(z,w)-  w\partial_{zz}^2 f_2(z,w)+4\pi i h z\partial_zf_2(z,w) +2\pi i h f_2(z,w) =0\,.
\end{equation}
which is equivalent through a change of variables~\amscite{PolyaninNazaikinskii16a}*{3.8.3.4}:
\begin{displaymath}
  (z,w,f_2)\rightarrow \left (\frac{z}{ w}, \frac{1}{w}, \frac{1}{\sqrt{w}} f_2\right) .
\end{displaymath}
to the free particle Schr\"o\-din\-ger equation.
The above discussion allows us to characterise the image space of the metamorphism:
\begin{thm}
  \label{th:metamorphism-image-characterisation}
  A function \(F(x,y,b,r)\) on \(\Space[+]{R}{4}\) is the metamorphism~\eqref{eq:metamorphism} of a function \(f\in \FSpace{L}{2}(\Space{R}{})\) if and if only  \(F\) satisfies to the following conditions:
  \begin{enumerate}
  \item
    \(F(x,y,b,r)\) is annihilated by operators \(\oper{C}_1\)~\eqref{eq:Cauchy-first}, \(\oper{C}_2\)~\eqref{eq:Cauchy-second} and \(\oper{S}_1\)~\eqref{eq:strctural-1}.
  \item
    \(F(\cdot,\cdot, b_0,r_0)\) is square-integrable in the sense of the inner product   \(\langle \cdot, \cdot \rangle_{(b_0, r_0)}\)~\eqref{scalar-product-b0-r0} for some \((b_0,r_0)\in\Space[+]{R}{2}\).
  \end{enumerate}
\end{thm}
\begin{proof}
  The necessity of both conditions was discussed above. For sufficiency, let a function \(G\) is annihilated by \(\oper{C}_1\), \(\oper{C}_2\) and \(\oper{S}_1\). If \(G(x,y, b_0,r_0) = 0\) for all \((x,y) \in \Space{R}{2}\) then the initial value problem for the parabolic equation~\eqref{eq:strctural-generic} implies \(G \equiv 0\) on \(\Space[+]{R}{4}\).

  Now, based on the square-integrability of \(F\) we use contravariant transform expression~\eqref{eq:contravariant-transform} to construct a function \(f \in \FSpace{L}{2}(\Space{R}{})\).
  Then \(G\coloneqq F - \metamorph{f}\) is annihilated by operators \(\oper{C}_1\), \(\oper{C}_2\), \(\oper{S}_1\) and \(G(x,y, b_0, r_0) = 0\) for all \((x,y) \in \Space{R}{2}\). Thus \(G\equiv 0\) (as explained above) and therefore \(F = \metamorph{f}\) on \(\Space[+]{R}{4}\).   
\end{proof}
Although the Gaussian and the metamorphism based on it are genuinely remarkable in many respects, other covariant transforms~\eqref{eq:covariant-tr-schrodinger-type} with fiducial vectors~\eqref{eq:fiducial-Er=0} and~\eqref{eq:fiducial-Er-not-0} deserve further attention as well.

\medskip

\section*{Acknowledgments}
\label{sec:acknowledgments}
Authors are grateful to Prof.~Alexey Bolsinov for useful comments and
suggestions on this work. Suggestions of anonymous referees helped to improve the presentation.\\
{\bf Funding and Conflicts of interests/Competing interests.}
The first named author was sponsored by the Albaha University
(SA).
The authors have no relevant financial or non-financial interests to disclose.

\small
\bibliography{abbrevmr,akisil,analyse,algebra,arare,aclifford,aphysics,acompute,ageometry,acombin}
\ifspringer
\else
\bibliographystyle{abbrv}
\fi

\end{document}